\documentclass[14pt]{article}
\usepackage[english]{babel}
\usepackage{amsfonts,amssymb,amsmath,amstext,amsbsy,amsopn,amscd,amsthm,graphicx,euscript,graphicx}
\usepackage{graphics}
\newtheorem{thm}{Theorem}
\newtheorem{lm}{Lemma}

\newtheorem{crl}{Corollary}

\newcounter{tdfn}
\newcounter{trk}
{\endtrivlist}

\def\:{\colon}

\def\R{{\mathbb R}}
\def\Z{{\mathbb Z}}
\def\0{{\mathbf 0}}
\def\1{{\mathbf 1}}

\def\C{{\mathbb C}}
\def\R{{\mathbb R}}

\title{Invariants of Classical Braids Valued in  $G_{n}^{2}$}

\author{V.O.Manturov \footnote{Research is carried out with the support of Russian Science Foundation (project no. 16-11-10291).}}

\date{}

 \def\R{{\mathbb R}}
 \def\Z{{\mathbb Z}}

\begin{document}
% \magstep2

\maketitle

AMS MSC 57M25, 57M27

{\Large

\begin{abstract}
The aim of the present note is to construct invariants of the Artin braid group valued in $G_{N}^{2}$, and further study of groups related to $G_{n}^{3}$. In the groups $G_{n}^{2}$, the word problem is solved; these groups are much simpler than $G_{n}^{3}$.
\end{abstract}

In \cite{Great}, the author introduced families of groups $G_{n}^{k}$ depending on two positive integers $n>k$ and formulated the following principle:
{\em if a dynamical system describing a motion of  $n$ particles admits a nice general position codimension $1$ property governed precisely by $k$ particles then this dynamical system has a topological invariant valued in $G_{n}^{k}$.}

In the paper \cite{MN}, the author calculated explicitly a partial case of this common principle: when we consider a continuous motion of $n$ pairwise distinct points on $\R^{2}$ and choose the property `` three points are collinear'', we get a homomorphism from the $n$ strand pure braid group $PB_{n}$ to the group $G_{n}^{3}$. This allows one to get powerful invariants of classical braids. Let us present the definition on  $`G_{n}^{3},n\ge 3$, which justifies the group $G_{n}^{3}$. The group $`G_{n}^{3}$ is given by a presentation having $3{n \choose 3}$ generators $a'_{ijk}$ indexed by triples of distinct numbers $i,j,k\in \{1,\cdots, n\}$ up to the reverse of order. Thus, $a'_{123}=a'_{321}\neq a'_{132}$. In the presentation
$$`G_{n}^{3}=\langle a'_{ijk}|(1),(2),(3)\rangle,$$
the relation (1) means $(a'_{ijk})^{2}=1$ for all pairwise distinct $i,j,k\in \{1,\cdots,n\}$; (2) means $a'_{ijk}a'_{pqr}=a'_{pqr}a'_{ijk}$ for all triples of pairwise distinct numbers $i,j,k$ and $p,q,r$ from $1$ to $n$, if $|\{i,j,k\}\cap \{p,q,r\}| <2$,
(3) means that for each quadruple of pairwise distinct indices $i,j,k,l$ from $1$ to $n$ the formula $(a'_{ijk}a'_{ijl}a'_{ikl}a'_{jkl})^{2}=1$ holds.

The group $G_{n}^{3}$ is obtained from $`G_{n}^{3}$ by identifying generators with coinciding triples of indices: the generator $a_{ijk}$ of $G_{n}^{3}$ equals $a'_{ijk}=a'_{jik}=a'_{ikj}$.

Note that $`G_{3}^{3}$ is the free product of three groups $\Z_{2}$, and the group $`G_{4}^{3}$ has no relations of type (2), since every two subsets of cardinality  3 of the set $\{1,2,3,4\}$ have intersection of cardinality at least $2$. When considering  $G_{n}^{3}$ one usually requires $n>3$; however, the group $`G_{3}^{3}$ is interesting as well as other groups $`G_{n}^{3},n>3$.

We shall also use the group $G_{n(n-1)}^{2}$; however, when defining this group (unlike the standard definitions of the group  $G_{N}^{2}$, see e.g., \cite{Bardakov,Great}) the indices $p,q$ of generators $a_{p,q}$ will be not elements of the set $\{1,\cdots, n(n-1)\}$ but rather pairs of distinct elements from $1$ to $n$. Hence, for generators of the group we can take, for instance, $a_{12,34},a_{12,31}$. The relations in this group are standard: $a_{p,q}^{2}=1$ for each generator $a_{p,q}$, $a_{p,q}a_{r,s}=a_{r,s}a_{p,q}$ for each pairwise distinct sets $p,q,r,s$ and $(a_{p,q}a_{p,r}a_{q,r})^{2}=1$ for pairwise distinct sets $p,q,r$.

Let us formulate the key lemma for our paper.
\begin{lm}
The map $\phi:`G_{n}^{3}\to G_{n(n-1)}^{2}$ given by the formula $a'_{ijk} \mapsto a_{ij,ik} a_{kj,ki},$  is well defined.\label{l1}
\end{lm}

First note that the two factors $a_{ij,ik}$ and $a_{kj,ki}$ commute, which means $\phi(a'_{ijk})=\phi(a'_{kji})$.

This lemma follows from a direct calculation. Let us check the most interesting relation:
$$(a'_{ijk}a'_{ijl}a'_{ikl}a'_{jkl})^{2}\mapsto (a_{ij,ik}a_{kj,ki}a_{ij,il}a_{lj,li}a_{ik,il}a_{lk,li}a_{jk,jl}a_{lk,lj})^{2}$$

$$=(a_{ij,ik}a_{ij,il}a_{ik,il})^{2}(a_{lk,lj}a_{lk,li}a_{lj,li})^{2}a_{jk,jl}^{2}a_{kj,ki}^{2}=1.$$

Note that this lemma is important by itself. For the group $G_{N}^{2}$, there is a simple minimality criterion for the word length (see, e.g., \cite{Coxeter}): a word $g$ in the standard presentation of $G_{N}^{2}$ has minimal length if and only if no word ${\tilde g}$ equivalent to $g$ by applying exchanges  $a_{p,q}a_{p,r}a_{q,r}\mapsto a_{q,r}a_{p,r}a_{p,q}$ and commutativity relations $a_{pq}a_{rs}\mapsto a_{rs}a_{pq}$ contains two adjacent identical letters $a_{p,q}a_{p,q}$ (in our case $N=n(n-1)$ and each letter $p,q,r,s$ itself consists of a pair of indices). Thus we get a  {\em sufficient minimality condition} for words from $`G_{n}^{3}$: if the image $\phi(\alpha)$ is minimal, then the word $\alpha$ is itself minimal.

The author does not know whether the map $\phi$ is an injection; this is an important open problem. It is also important to know whether one can construct some map analogous to the map $\phi$, taking $G_{n}^{3}$ to some group $G_{M(n)}^{2}$ (for some large $M(n)$ depending on $n$). A positive answer to these questions could shed light on the solvability problem in the groups $G_{n}^{k}$ for $k>2$.

Now let us construct the map  $f$, which associates with a pure braid $\beta\in PB_{n}$ ($n\ge 3$) an element from $`G_{n}^{3}$. We shall deal with braids for which the initial and the final set of point are uniformely distributed along the circle: $z_{j}=exp(\frac{2\pi j}{n})$. By a  {\em braid} we mean a set of smooth functions  $\beta(t)=\{z_{1}(t),\cdots,z_{n}(t)\}, t\in [0,1]$ valued in $\C^{1}=\R^{2}$ such that $b(0)=b(1)$ coincides with the set give above, so that all $z_{i}(t)$ are pairwise distinct for all $t$. A {\em critical moment} is such a value $t$, for which there exist some three indices $i,j,k$ such that $z_{i}(t),z_{j}(t),z_{k}(t)$ are collinear.
We say that a braid is {\em good and stable} if for this braid:
\begin{enumerate}
\item the set of critical moments is finite;
\item в for each critical moment $t$ there exist exactly one triple of indices $(i,j,k)$ for which $z_{i}(t),z_{j}(t),z_{k}(t)$ are collinear;
\item any small perturbation of the braid does not change the number of critical moments.
\end{enumerate}

By a small perturbation we can make any braid good and stable.

With a good and stable pure braid $\beta$ one naturally associates a word in generators $a'_{ijk}$: with each critical moment $t_{l}$ with corresponding three collinear points $(i_{l},j_{l},k_{l})$ with  $j_{l}$ in the middle, we associate the generator $a'_{i_{l},j_{l},k_{l}}$. The word $f(\beta)$ is the product of all generators corresponding to all critical moments (counted as $t$ grows from $0$ to $1$).

\begin{thm}
The map constructed above $\beta\mapsto f(\beta)$ is a homomorphism from the pure braid group $PB_{n}$ to the group $`G_{n}^{3}$.
\end{thm}
The proof of this theorem essentially repeats the proof of the theorem from \cite{MN} about the map from $PB_{n}$ to $G_{n}^{3}$. We shall use the standard principle from \cite{Great} saying that for the isotopy of two braids it suffices to consider only singularities of codimension two. As singularities of codimension one (triples of collinear points) give rise to generators, singularities of codimension two lead to relations. The latter look as follows (see \cite{MN}):

\begin{enumerate}
\item a nonstable triple point, which disappears after a small perturbation; this leads to ${a'_{ijk}}^{2}=1$;

\item coincidence of two moments when two independent triple points appear. This corresponds to 
$a'_{p}a'_{q}=a'_{q}a'_{p}$, where two triples $p,q$ have at most one common index;

\item four collinear points; this yields the relation when on the LHS we have a product of four generators and on the RHS we have the product of the same generators in the inverse order. Another way to describe this relation is the cyclical deformation of the dynamics containing of eight elementary deformations.
    
    These four generators correspond to all possible triples of indices which are subsets of the given four indices; for example, $a'_{ijk},a'_{ijl},a'_{ikl},a'_{jkl}$. The only novelty in comparison with \cite{MN} is that the relation in the group $`G_{n}^{3}$ is more precise  than that in $G_{n}^{3}$. If we consider the line passing through the four points and order denote these numbers of these points as they are along the line $i,j,k,l$, we can check directly that in the above cyclic deformation the generator  $a'_{ijk}$ (corresponding to the three points on the line with $j$ in the middle) can not be adjacent with the generator $a'_{ikl}$. Hence, we get a certain  word which can be obtained from the LHS of the relation (3) by some cyclic perpmutation and order reversal. Since the squares of all generators are equal to $1$, the resulting word obtained by order reversal is the inverse element of the group.
\end{enumerate}

The product of pure braids leads to the product of words.

\begin{crl}
The map $\Phi=\phi\circ f$ is a homomorphism $PB_{n}\to G_{n(n-1)}^{2}$.
\end{crl}
As we see from the construction above, the map $\Phi$ can be constructed without mentioning the ``intermediate'' group $`G_{n}^{3}$: for a general position dynamics one can directly write down an element of $G_{n(n-1)}^{2}$.

The mapping $f$ constructed above is a generalization of rather a strong invariant of classical braids constructed inn \cite{MN}. The composition $\Phi=\phi\circ f$ for braids can be readily calculated. Say, with a generator of the pure three-strand braid group where the first and the second points stay fixed and the third one goes around the second one, one associates the product of four generators $G_{n(n-1)}^{2}$: the two moments corresponding to triples of collinear points in different orders, lead to a word of length  $4$ in the group $`G_{3}^{3}=\Z_{2}*\Z_{2}*\Z_{2}$. This word is minimal.

We conclude our note by constructing a homomorphism from $`G_{n}^{3}$ to the automorphism group of the free product of several groups $\Z_{2}$, which is ``spiritually'' similar to the Hurwitz action of the Artin braid group on the free group by automorphisms.

Namely, we define the image  $g=g(a_{ijk})$ as $g: a_{ij}\mapsto a_{ik}a_{ij}a_{ik},a_{kj}\mapsto a_{ki}a_{kj}a_{ki},
a_{m}\mapsto a_{m}, \mbox{where } m\neq \{ij\},\{ki\}$

\begin{thm}
The map  $g$ is a well defined homomorphism.
\end{thm}

The statement follows from a direct check of the relations of the group $`G_{n}^{3}$. One can mention that the check of the main relation (3) is similar to the proof of Lemma \ref{l1} and the third Reidemeister move.

}

 \end{document}